\documentclass[12pt, a4paper]{amsart}
\usepackage[hmargin=30mm, vmargin=25mm, includefoot, twoside]{geometry}
\usepackage[bookmarksopen=true]{hyperref}
\usepackage[latin1]{inputenc}  % accents 8 bits dans la source
\usepackage[french,english]{babel}
\usepackage{amssymb,verbatim,  pstricks, pst-plot, textcomp, fontenc, euscript}
\usepackage[alphabetic]{amsrefs}
\usepackage{latexsym}
\usepackage{mathrsfs}
\usepackage{xspace}
\usepackage{enumerate, paralist}
\newtheorem{prop}{Proposition}%[section]
\newtheorem{thm}[prop]{Theorem}
\newtheorem*{thm*}{Theorem}
\newtheorem*{prop*}{Proposition}

\newtheorem*{addendum*}{Addendum}

\newtheorem*{convention*}{Convention}
\theoremstyle{definition}
\newtheorem*{defn*}{Definition}
\newtheorem{defn}[prop]{Definition}
\newtheorem{remark}[prop]{Remark}

\newtheorem{question}[prop]{Question}
\newtheorem*{scholium*}{Scholium}
\theoremstyle{remark}
\newtheorem{example}[prop]{Example}
\newtheorem*{example*}{Example}
\numberwithin{equation}{section}

\newcommand{\FF}{\mathbf{F}}

\newcommand{\cat}{{\upshape CAT($0$)}\xspace}

\newcommand{\tangle}[2]% angle de Tits
{\angle_\mathrm{T}(#1,#2)}
\newcommand{\aangle}[3]% angle d'Alexandrov
{\angle_{#1}(#2,#3)}
\newcommand{\cangle}[3]% angle de comparaison
{\overline{\angle}_{#1}(#2,#3)}
\begin{document}

\title[Simplicity of some twin tree lattices]{Simplicity of some twin tree automorphism groups\\ with trivial commutation relations}

\author[J.~Morita]{Jun Morita*}
\address{Institute of Mathematics, University of Tsukuba, 1-1-1 Tennodai, Tsukuba, Ibaraki 305-8571 Japan}
\email{morita@math.tsukuba.ac.jp}
\thanks{*Supported in part by Grant-in-aid for Science Research (Monkasho Kakenhi) in Japan}

\author[B. R\'emy]{Bertrand R\'emy**}
\address{Universit\' e Lyon 1\\
Institut Camille Jordan\\
UMR 5208 du CNRS\\
43 blvd du 11 novembre 1918\\
F-69622 Villeurbanne Cedex, France}
\email{remy@math.univ-lyon1.fr}
\thanks{**Supported in part by Institut Universitaire de France}

\date{\today}

\maketitle

\begin{abstract}
We prove simplicity for incomplete rank 2 Kac-Moody groups over algebraic closures of finite fields with trivial commutation relations between root groups corresponding to prenilpotent pairs.
We don't use the (yet unknown) simplicity of the corresponding finitely generated groups (i.e., when the ground field is finite).
Nevertheless we use the fact that the latter groups are just infinite (modulo center).
\end{abstract}

\bigskip

\section*{Introduction}
\label{s - intro}

In this paper we prove simplicity (up to center) of some (incomplete) Kac-Moody groups over algebraic closures of finite fields.
At first glance, this might be a surprising result because the examples which are usually given to introduce incomplete Kac-Moody groups (as defined by J.~Tits \cite{TitsKM}) are of affine type, and the latter groups have a matrix interpretation.
For instance, a Kac-Moody group of type $\widetilde{{\rm A}}_n$ over some field $F$ is isogenous to ${\rm SL}_n(F[t,t^{-1}])$.
In fact, any $F$-split simple algebraic group ${\bf G}$ gives rise to a Kac-Moody group functor $R \mapsto {\bf G}(R[t,t^{-1}])$ on $F$-algebras.
The values over fields of such a functor are (Kac-Moody) groups admitting a lot of (congruence) quotients, since the ring $R[t,t^{-1}]$ has arbitrarily small ideals.

\smallskip 

The question is thus: given a certain class of ground fields, which types of Kac-Moody groups shall we exclude to hope for simplicity?
The situation over finite ground fields is almost completely understood \cite{CaRe}.
The outcome suggests that among the irreducible generalized Cartan matrices, the only types that should be excluded are the affine ones.
To be more precise, this general picture over finite fields is completely confirmed except when the generalized Cartan matrix is $2 \times 2$, in which case the problem is only half solved \cite{CaReRk2}.
The connection with our case, where ground fields are of the form $\overline{{\bf F}_q}$, is that simplicity over finite ground fields easily implies simplicity over the algebraic closure (\ref{ss - proof}, Remark \ref{rk - simple}).
We deal here with the only case where simplicity over finite ground fields is still an open question.

\begin{thm*}
Let $A = \left( \begin{array}{cc} \hfill 2  & -n \\ -m & \hfill 2 \end{array} \right)$ be a generalized Cartan matrix of indefinite type, i.e. $mn>4$.
Let $\mathscr{G}_A$ be the corresponding simply connected incomplete Kac-Moody group functor and let $F$ be an algebraic closure of a finite field.
Assume that $m,n \geqslant 2$.
Then, the group $\mathscr{G}_A(F)/Z\bigl( \mathscr{G}_A(F) \bigr)$ is simple.
\end{thm*}

In particular, this theorem settles the last case needed to prove the following statement (see Remark \ref{rk - simple final}, subsection \ref{ss - heuristic}): 
{\it irreducible, simply connected, non-affine Kac-Moody groups over algebraic closures of finite fields are simple modulo their centers}.
The picture over finite fields is slightly less complete.

\smallskip

The reason why excluding affine types for simplicity over finite ground fields has a geometric explanation which naturally leads us to introduce the main tool in the investigation of these groups, namely buildings (another concept introduced by J.~Tits and presented in \ref{ss - twin buildings}). 
Roughly speaking a building is a nice, symmetric, simplicial complex designed to admit group actions.
By definition, a building is covered by subcomplexes (called apartments) which are all isomorphic and whose geometry is fully encoded by a Coxeter group which is called the Weyl group of the building.
An infinite Weyl group is a Euclidean reflection group if and only if it has polynomial growth for its natural generating set.
For generalized Cartan matrices of size $\geqslant 3$, simplicity occurs (at least over finite fields) precisely when the Weyl group of the buildings is not Euclidean, because then the associated root system has some nice weak hyperbolicity properties.
The proof of our result is also related to some kind of hyperbolicity since our assumption $mn>4$ corresponds to hyperbolic root systems of rank 2.
This proof requires in addition the use of some weak version of simplicity, called the normal subgroup property, which is reminiscent to a famous result of G.~Margulis about lattices in higher rank Lie groups.

\smallskip

The structure of the paper is the following.
In Section 1, we introduce the basic objects used in the paper, namely twin buildings and Kac-Moody groups.
In Section 2, we recall the situation over finite fields because we need to state the normal subgroup property in this case.
In Section 3, we prove our main theorem and mention the remaining related problem for finitely generated Kac-Moody groups.

\smallskip

Let us finally introduce some notation. 
Concerning groups, $Z(G)$ means the center of a group $G$.
Concerning rings, ${\bf Z}$ (resp. ${\bf Q}$, ${\bf R}$) means the set of integral (resp. rational, real) numbers.
In this article, $p$ is a prime number and $q$ a power of some $p$; at last, ${\bf Q}_p$ (resp. ${\bf F}_p$, ${\bf F}_q$) means the field of $p$-adic numbers (resp. a prime field of characteristic $p$, a finite field of order $q$).

\bigskip

\section{Twin building and Kac-Moody theory}
\label{s - TB & KM}
All the theories in this subsection are due to J.~Tits, see for instance \cite{TitsVancouver} and \cite{TitsTwin} for (twin) buildings and \cite{TitsKM} for Kac-Moody groups.

\subsection{Twin building theory}
\label{ss - twin buildings}

Let us first recall the definition of a building.
If $W = \langle s \in S \mid (st)^{M_{st}}=1\rangle$ is a Coxeter group defined by the Coxeter matrix $[M_{st}]_{s,t \in S}$, there is a simplicial complex, called the {\it Coxeter complex}~$\Sigma$ of $(W,S)$, on the maximal simplices of which $W$ acts simply transitively \cite{BBK}. 
In this context, simplices are rather called {\it facets}.
Coxeter complexes (seen as simplicial or metric spaces) are generalized tilings on which the initial Coxeter group acts as a generalized reflection group generated by natural involutions (reflections in faces of a given {\it chamber}, i.e. a maximal facet).

\smallskip

Up to removing the facets with infinite stabilizers, there exists a geometric realization for $\Sigma$, usually different from the one introduced in Bourbaki, carrying a complete metric such that the resulting metric space is non-positively curved and contractible.
Technically the notion is that of a complete {\it \cat-space}.
Since we will use this terminology without going into technical details, we simply refer to \cite{BH}.

\begin{defn}
A {\it building of type $\Sigma$}~is a simplicial complex covered by sub-complexes all isomorphic to the Coxeter complex $\Sigma$, called {\it apartments} and required to satisfy the following axioms.
\begin{enumerate}
\item[(i)]~Any two simplices are always contained in an apartment.
\item[(ii)]~Given any two apartments $A$, $A'$ there is an isomorphism $A \simeq A'$ fixing $A \cap A'$.
\end{enumerate}
The group $W$ is called the {\it Weyl group}~of the building.
\end{defn}

The above axioms can be motivated by metric considerations.
Indeed, they can be used to glue together the above (Davis-Moussong) metrics on each apartment in order to define a complete \cat metric on the building: axiom (i) says that computing the distance between two points can always be done by doing it in a suitable apartment and axiom (ii), up to additional work in order to define suitable 1-lipschitz retractions, shows that the distance computed this way doesn't depend on the apartment.

\begin{example}
Let $D_\infty$ be the infinite dihedral group, i.e. the group generated by two reflections in consecutive integers on the real line.
Then a building of type $D_\infty$ is a tree (without pending leaf).
Note that such a tree may have no automorphism at all since trees in which any two vertices have distinct valencies are not excluded by the axioms (the isomorphism in (ii) need not be  defined globally).
\end{example}

\begin{example}
The Coxeter complex of type $\widetilde A_2$ is the one given by the tiling of ${\bf R}^2$ by regular triangles.
Buildings whose apartments have this shape are called triangle buildings; they appear as Bruhat-Tits buildings for Lie groups like ${\rm SL}_3$ over local fields.
More generally, one consequence of Cartan and Bruhat-Tits theories is the possibility to associate to any $S$-arithmetic group a complete \cat-space on which it acts nicely.
These spaces are obtained as products of symmetric spaces and of {\it Euclidean buildings}, i.e. buildings in which apartments are Euclidean tilings.
\end{example}

New interesting questions occur when the buildings of the geometric actions under consideration are no longer Euclidean.
Many examples of buildings with hyperbolic tilings as apartments are available thanks to Kac-Moody theory (\ref{ss - KM}).
Such exotic buildings provide more opportunities to study non-linear discrete groups via geometric actions.

\smallskip

Let us now turn very quickly to twinnings \cite{TitsTwin}.
Initially, the idea is to extend some rigidity properties (useful in the classification of spherical buildings)  to non-spherical buildings, provided they are twinned with another isomorphic building.
The idea to add a twin to a non-spherical building enables one to define an opposition relation between facets in the two different buildings.
Conventionally, each of the two twinned buildings is given a sign $\pm$.
This opposition relation between facets of opposite signs is a substitute for  the existence of a longest element in the Weyl group of the buildings.

\begin{example}
By Bruhat-Tits theory, the groups ${\rm SL}_n\bigl( \FF_q (\!( t^{\pm 1} )\!) \bigr)$, where $\FF_q (\!( t^{\pm 1} )\!)$ are locally compact non-Archimedean fields of formal Laurent series, act on isomorphic Euclidean buildings, say $X_\pm$. 
There is a natural twinning between $X_-$ and $X_+$ such that the discrete group ${\rm SL}_n(\FF_q[t,t^{-1}])$ (embedded diagonally in the product of the two previous groups) acts on $X_- \times X+$ and preserves opposition of chambers of opposite signs.
\end{example}

Given a homogeneous tree, there are uncountably many ways to associate to it a twin tree, but most twinnings have no automorphism at all \cite{RonanTitsIJM}.
Still, the additional Moufang condition on twin buildings guarantees the existence of enough automorphisms for these buildings.
We will not go into details, but we simply mention that Kac-Moody theory provides lots of examples of twin buildings satisfying the Moufang condition.
Even more exotic (i.e., non Kac-Moody) Moufang twin buildings with enough automorphisms are also available by means of more down-to-earth constructions, see \cite{RemRon} and also \cite{AbRe}.

\subsection{Kac-Moody groups}
\label{ss - KM}
Kac-Moody groups are constructed from the same kind of data as for Chevalley groups, namely a ground field and some Lie-theoretic data classifying semisimple Lie algebras \cite{RemAst}.

\smallskip

More precisely, a {\it generalized Cartan matrix}~is an integral matrix $A = [A_{s,t}]_{s,t \in S}$ indexed by a set $S$ (which is here assumed to be finite), such that $A_{s,s}=2$ for any $s \in S$ and $A_{s,t} \leqslant 0$ for any $s \neq t$ in $S$; it is further required that $A_{s,t} = 0$ if and only if $A_{t,s} =0$.
From this Lie-theoretic matrix, a certain group functor over rings can be constructed by generators and relations \cite{TitsKM}.
It is a heavy machinery of algebraic and combinatorial nature, which gives a Chevalley group scheme if the matrix $[A_{s,t}]_{s,t \in S}$ is a Cartan matrix (i.e., if it can be written as the product of a diagonal matrix with a positive definite symmetric matrix).
In fact, as in this classical case, the matrix $A$ only determines Kac-Moody group functors up to center, and in what follows we always use the simply connected groups (this choice plays no significant role for our purposes -- it makes simplicity results easier to state).
A {\it Kac-Moody group}~is the value of a Kac-Moody functor on a field, called the {\it ground field}~in what follows.

We are interested in the geometric outcome of this construction: {\it a Kac-Moody group acts on the product of two Moufang twin buildings and the kernel of the action is its center}.
It is a well-known fact that a group enjoying the structure of a Tits system (also called BN-pair) naturally acts (strongly transitively) on a building.
In the case of a Kac-Moody group of non-classical (i.e., non-Chevalley) type, there are two conjugacy classes of subgroups which leads to two distinct buildings.
Moreover the Weyl group, i.e. the shape of the apartments of the twinned buildings $X_\pm$, is explicitly known since its Coxeter matrix $[M_{s,t}]_{s,t \in S}$ is determined by the rule $M_{s,t} = 2$ (resp. $3,4,6$ or $\infty$) according to whether $A_{s,t}A_{t,s}$ is equal to $0$ (resp. $1,2,3$ or is $\geqslant 4$).
At last the associated buildings $X_\pm$ are locally finite if and only if the ground field is finite, which we assume until the end of the next section.
This implies that, for the \cat-metric, the isometry groups ${\rm Isom}(X_\pm)$ are locally compact for the compact open topology, and as such admit Haar measures.

\begin{example}
Over a given field ${\bf F}$, for a suitable choice of generalized Cartan matrices (namely for those of affine type), the corresponding Kac-Moody groups are of the form ${\bf G}({\bf F}[t,t^{-1}])$ where ${\bf G}$ is a semisimple algebraic group over ${\bf F}$.
Then the actions of ${\bf G}({\bf F}[t,t^{-1}])$ on the associated twin buildings are those given by Bruhat-Tits theory by seeing ${\bf G}({\bf F}[t,t^{-1}])$ as a subgroup of the two completions 
${\bf G}\bigl({\bf F}(\!(t)\!)\bigr)$ and ${\bf G}\bigl({\bf F}(\!(t^{-1})\!)\bigr)$.
\end{example}

\begin{example}
Using the rule $[A_{s,t}]_{s,t \in S} \to [M_{s,t}]_{s,t \in S}$, we easily see that many buildings whose apartments are real hyperbolic tilings are made available by Kac-Moody theory.
An interesting point in this construction is the fact that a Kac-Moody group acts on each of the two twinned buildings in a highly transitive way (in particular it acts on each factor with a chamber as fundamental domain).
\end{example}

\bigskip

\section{Finitely generated Kac-Moody groups}
\label{s - fg KM}

In this section, we recall the general situation of twin building lattices, as far as the question of simplicity is concerned.
For this we need to recall some general notions from geometric group theory.

\subsection{A glimpse of geometric group theory}
\label{ss - GGT}
Roughly speaking, arithmetic groups are matrix groups with coefficients in rings of integers of global fields and in natural generalizations; examples of such groups are ${\rm SL}_n({\bf Z})$ or ${\rm SL}_n({\bf Z}[1/p])$.
An arithmetic group appears as a subgroup in a product of (real and totally disconnected) Lie groups (e.g., ${\rm SL}_n({\bf Z}) < {\rm SL}_n({\bf R})$, and ${\rm SL}_n({\bf Z}[1/p]) < {\rm SL}_n({\bf R}) \times {\rm SL}_n({\bf Q}_p)$ for the diagonal inclusion).
Furthermore a non-compact simple Lie group naturally acts on a complete \cat-space \cite{SMF18}.
It is a symmetric space if the simple Lie group is defined over the real numbers.
When the ground field of the simple Lie group is a non-Archimedean local field, the metric space is a Euclidean building (the construction of the latter space is not trivial at all, it follows from the so-called Bruhat-Tits theory \cite{Rousseau}).
Putting these two facts together (and forgetting the step involving the ambient topological groups), we obtain an interesting situation (called here a {\it geometric action}) where a discrete group $\Gamma$ acts on a metric space $(X,d)$ so that:

\smallskip

\begin{enumerate}
\item[(GA1)] the metric $d$ on $X$ is complete and \cat;
\item[(GA2)] the group $\Gamma$ acts by isometries and properly discontinuously on $X$;
\item[(GA3)] the $\Gamma$-action has a nice fundamental domain.
\end{enumerate}

\noindent By "nice" fundamental domain, we can mean compact, but compactness is usually too strong.
More technically, it means that the full isometry group ${\rm Isom}(X,d)$ carries a Haar measure and that the corresponding invariant measure on the homogeneous space ${\rm Isom}(X,d)/\Gamma$ has finite volume.
We say then that $\Gamma$ is a {\it lattice}~for $(X,d)$.

\begin{example}
\label{ex - Poincare}
The symmetric space associated to ${\rm SL}_2({\bf R})$ is Poincar\'e's upper half-plane $\mathbb{H}^2_{\bf R}$ and the group ${\rm SL}_2({\bf Z})$ acts on it with the well-known fundamental domain $\{ z \in {\bf C} \,\, : \,\, \mid\! z \!\mid \,\geqslant 1$ and $\mid\! {\rm Re}(z) \!\mid \,\leqslant {1 \over 2}\}$.
\end{example}

\begin{example}
The Bruhat-Tits building associated to the rank 1 non-Archimedean simple Lie group ${\rm SL}_2({\bf Q}_p)$ is the homogeneous tree $T_{p+1}$ of valency $p+1$.
The natural action of the lattice ${\rm SL}_2({\bf Z}[{1 \over p}])$ is the diagonal action on the mixed product $\mathbb{H}^2_{\bf R} \times T_{p+1}$ of a differentiable manifold and a simplicial complex.
\end{example}

\begin{example}
To obtain a geometric action of a lattice on a product of two trees, one can use slightly less familiar matrix groups.
Namely, start with a quaternion algebra over ${\bf Q}$, say $H$, such that $H({\bf R})$ is a skew-field (in arithmetic terms, $H$ is ramified at $\infty$); pick two prime numbers $p$ and $l$ such that $H({\bf Q}_p)$ and $H({\bf Q}_l)$ are matrix algebras.
Then the elements in $H\bigl({\bf Z}[{1\over pl}]\bigr)$ form a discrete group having a geometric action on $T_{p+1} \times T_{l+1}$, and the fundamental domain is compact.
\end{example}

A typical question in geometric group theory consists in asking what can be said about a discrete group once it is known to admit a geometric action on a particularly nice \cat-space (e.g. a non-spherical building or a cube complex -- products of trees belong to both classes).
Relevant questions are for instance related to freeness, linearity, residual finiteness, simplicity etc.
The historical statement, in connection with Example \ref{ex - Poincare}, is the proof that ${\rm SL}_2({\bf Z})$ contains a finite index subgroup isomorphic to the free group $F_2$ (this is F.~Klein's ping-pong argument).

\subsection{Non-affine higher-rank finitely generated Kac-Moody groups}
\label{ss - NSP}
We can now go back to the objects defined in the previous section.
Let $\Lambda$ be a Kac-Moody group over a finite field $\FF_q$ of order $q$.
Then the diagonal $\Lambda$-action on $X_-\times X_+$ is geometric in the sense of the axioms (GA) in \ref{ss - GGT}.
Using the group combinatorics of twinTits systems, we can see that a fundamental domain is given for instance by the product of a negative chamber by a suitable positive apartment.
The starting point of the analogy between Kac-Moody groups over finite fields and $S$-arithmetic groups is the following result \cite{RemCRAS}: {\it at least when $q > \#S$, the group $\Lambda$, which is finitely generated by construction, is a lattice in ${\rm Isom}(X_\pm) \times {\rm Isom}(X_\pm)$}.
In fact, the covolume of $\Lambda$ is given by $\sum_{w \in W} q^{-\ell(w)}$ for a suitable normalization of Haar measures; in particular, for twin trees (where $\# S=2$) the covolume is always finite since $W$ has linear growth in that case.

Now the main structure result on normal subgroups of lattices in Lie groups is due to G.~Margulis \cite[Lecture 4]{Benoist}:
{\it let $\Gamma$ be an irreducible lattice in a higher-rank semisimple Lie group.
Then for any $\Delta \triangleleft \Gamma$, either the subgroup $\Delta$ is finite and central, or $\Delta$ has finite index in $\Gamma$.}
A group all of whose normal subgroups satisfy the previous dichotomy is said to have the {\it normal subgroup property}, (NSP) for short.
This is a typical result to try to generalize for lattices in products of buildings obtained from Moufang twin buildings.
This was indeed checked in \cite{RemInt}: {\it let $\Lambda$ be an irreducible Kac-Moody group over a finite field.
Then $\Lambda$ has {\rm (NSP)}~whenever it is a lattice of the product of its two twinned buildings}~(i.e., whenever the finite ground field is big enough with respect to the growth of the Weyl group -- see the above covolume formula).
The proof follows Margulis' general strategy consisting in showing that for an infinite $\Delta \triangleleft \Lambda$, the discrete group $\Lambda/\Delta$ is both amenable and Kazhdan (implying compactness, hence finiteness by discreteness).

The next step after (NSP) is simplicity.
Here is a simplified statement of what is proved in \cite{CaRe}: {\it let $\Lambda$ be a (simply connected) Kac-Moody group over the finite field $\FF_q$.
Assume that the generalized Cartan matrix defining $\Lambda$ is non-affine and indecomposable, say of size $n$. 
Then $\Lambda/Z(\Lambda)$ is simple whenever $q>n>2$.}
For this simplicity theorem, by (NSP) the key point is to rule out also the possibility to have finite quotients either; this is where the new conditions on the generalized Cartan matrix appear (non-affineness and $n>2$).
Indeed the argument to exclude finite quotients for $\Lambda$ uses the geometry of the root system of the Weyl group, more precisely the fact that whenever an infinite Coxeter group is irreducible, non-affine and of rank $>2$, then its root system has many hyperbolic triples:~seeing roots as half-spaces bordered by fixed-point sets of reflections in the Coxeter complex $\Sigma$, this means existence triples of pairwise disjoint roots in $\Sigma$ (which is clearly excluded for Euclidean reflection groups).

\begin{remark}
It is interesting to have simple groups occurring as lattices in products of buildings in which some freedom for the shape of the apartments is available.
Indeed, this leads to the following statement in geometric group theory \cite{CaReQI}: {\it there exists infinitely many quasi-isometry classes of finitely presented simple groups}.
\end{remark}

\bigskip

\section{Simplicity for non locally finite twin trees}
\label{s - simplicity}

We prove simplicity for hyperbolic rank 2 Kac-Moody groups over algebraic closures of finite fields (cf. Theorem in Introduction). 
This can be easily established when the corresponding Kac-Moody group is simple over a finite subfield (see Remark \ref{rk - simple}), so we concentrate on the case where the latter simplicity is still unknown.
This is when the commutation relations between root groups indexed by prenilpotent pairs are trivial.

\subsection{Simplicity without using simplicity}
\label{ss - proof}
Let us recall Tits functors, $\mathscr{G}_A$ as well as $\mathscr{T}_A$, associated with generalized Cartan matrices $A$
to produce the corresponding Kac-Moody groups \cite{TitsKM}.

\smallskip

Let $A = \left( \begin{array}{cc} \hfill 2  & -n \\ -m & \hfill 2 \end{array} \right)$ be a generalized Cartan matrix of indefinite type (i.e. $mn>4$).
We assume that $m,n \geqslant 2$, which implies that the commutation relations between root groups indexed by prenilpotent pairs are trivial \cite{Morita}.
Then, we obtain the corresponding Kac-Moody group $\mathscr{G}_A(F)$ over $F = \overline{{\bf F}_q}$ and the so-called standard maximal split torus
$\mathscr{T}_A(F) \simeq {\rm Hom}_{\mathbb{Z}}(\mathbb{Z}^2,F^\times)$.
The group $\mathscr{G}_A(F)$ is generated by root subgroups $U_\delta$ for all real roots $\delta$ in this case.

\smallskip

For each real root $\delta$, there is a natural isomorphism from the additive group $(F,+)$ onto $U_\delta$, which we denote by $r \mapsto u_\delta(r)$.
Tits' presentation \cite{TitsKM} implies that the group $S_\delta = \langle U_\delta, U_{-\delta} \rangle$ is isomorphic to ${\rm SL}_2(F)$ via an isomorphism
which sends $u_\delta(r)$ (resp. $u_{-\delta}(r)$) to 
$\left( \begin{array}{cc} \hfill 1  & r \\ 0 & \hfill 1 \end{array} \right)$ (resp. $\left( \begin{array}{cc} \hfill 1  & 0 \\ r & \hfill 1 \end{array} \right)$).
Then, $\mathscr{T}_A(F)$ is generated by $h_\delta(\mu)$ for all real roots $\delta$ and for all $\mu \in F^\times$, where $h_\delta(\mu)$ is an element of $S_\delta$ corresponding
to $\left( \begin{array}{cc} \hfill \mu  & 0 \\ 0 & \hfill \mu^{-1}  \end{array} \right)$.

Let $\alpha$ and $\beta$ be the simple roots defined by this presentation.
For any nonzero $j \in \mathbf{Z}$ we set $\gamma_j = \tau^j . \alpha$ where $\tau = w_\alpha(1)w_\beta(1)$ and $w_\delta(1) = u_\delta(1)u_{-\delta}(-1)u_\delta(1)$;  there exist integers $a_j$ and $b_j$ with $a_j b_j > 0$ such that 

\medskip

\centerline{$\gamma_j = a_j \alpha + b_j \beta$.}

\medskip

Note that in the geometric realization of the Weyl group $D_\infty$, the Coxeter complex (hence any apartment) is the real line. 
The reflections in the Weyl group are those with respect to the integers and the element $\tau$ acts as a translation along this line.

\smallskip

A element $t \in \mathscr{T}_A(F)$ given by this presentation has the form $t = h_\alpha(\mu) h_\beta(\nu)$, where $\mu, \nu \in F^\times$ are two multiplicative parameters.
Then, we have: 

\medskip

\centerline{$\alpha(t) = \mu^2 \nu^{\alpha(\beta^\vee)} = \mu^2 \nu^{-m}$}

\medskip

\noindent
and

\medskip

\centerline{$\gamma_j(t) = \mu^{\gamma_j(\alpha^\vee)} \nu^{\gamma_j(\beta^\vee)} 
= \mu^{2a_j+\beta(\alpha^\vee)b_j}\nu^{\alpha(\beta^\vee)a_j + 2b_j}
= \mu^{2a_j-nb_j}\nu^{-ma_j+2b_j}$,}

\medskip

\noindent
where $\gamma^\vee$ denotes the coroot of a real root $\gamma$.

\bigskip

{\it Proof of the theorem.}~
Let $K \triangleleft \mathscr{G}_A(F)$ be a non-central normal subgroup.
In order to prove our simplicity theorem (see Introduction), we must show that we have in fact  $K = \mathscr{G}_A(F)$.

\smallskip

Since each root  subgroup is conjugate to $U_\alpha$ or $U_\beta$, and since $S_\alpha \simeq S_\beta \simeq {\rm SL}_2(F)$, it is enough to show that $U_\alpha \cap K \neq \{ 1 \}$ and $U_\beta \cap K \neq \{ 1 \}$ (the group ${\rm SL}_2(F)$ doesn't contain any proper normal subgroup intersecting non-trivially a root group).

\smallskip

Since $F = \bigcup_{i \geqslant 1} {\bf F}_{q^i}$, we have $Z\bigl( \mathscr{G}_A(F) \bigr) = \bigcup_{i \geqslant 1} Z\bigl( \mathscr{G}_A({\bf F}_{q^i}) \bigr)$, and therefore 
there exists $\ell \geqslant 1$ such that ${\bf F}_{q^\ell} \subset F$ and $K \cap \mathscr{G}_A({\bf F}_{q^\ell})$ is non-central.
By the normal subgroup property \cite{RemInt}, and assuming that $\ell$ is large enough, 
the normal subgroup $K \cap \mathscr{G}_A({\bf F}_{q^\ell})$ has finite index, say $k$, in $\mathscr{G}_A({\bf F}_{q^\ell})$.
This implies, in particular, that $[\langle \tau \rangle : K \cap \langle \tau \rangle]$ divides $k$, which follows from

\smallskip

\centerline{$\displaystyle{
\begin{array}{lll}
k & = & [ \mathscr{G}({\bf F}_{q^\ell}) : K \cap \mathscr{G}({\bf F}_{q^\ell}) ]\\
& = &
[ \mathscr{G}({\bf F}_{q^\ell}) : \langle \tau \rangle (K \cap \mathscr{G}({\bf F}_{q^\ell}) ) ] \times
[ \langle \tau \rangle (K \cap \mathscr{G}({\bf F}_{q^\ell}) ) : K \cap \mathscr{G}({\bf F}_{q^\ell}) ]
\end{array}
}$}

\smallskip

\noindent
and

\smallskip

\centerline{$
[ \langle \tau \rangle (K \cap \mathscr{G}({\bf F}_{q^\ell}) ) : K \cap \mathscr{G}({\bf F}_{q^\ell}) ]
= [\langle \tau \rangle : K \cap \langle \tau \rangle]
$,}

\smallskip

\noindent
so that $\tau^k \in K$.
As a consequence, we have $[\tau^k, U_\alpha] \subset K$.

Let us start with $u \in U_\alpha-\{ 1 \}$, i.e. $u = u_\alpha(c)$ for some $c \in F^\times$.
It follows from the defining relations of an incomplete Kac-Moody group that we have $\tau^j U_\delta \tau^{-j} = U_{\tau^j.\alpha}$, so that: 

\medskip

\centerline{$[\tau^j,u] = (\tau^j u_\alpha(c) \tau^{-j}) u_\alpha(-c) = u_{\tau^j.\alpha}(r) u_\alpha(s)$}

\medskip

\noindent
for some suitable $r, s \in F^\times$.
Hence we see that for suitable powers $j$ (e.g. $j$ divisible by $k$) we can find elements in $\bigl( (U_\alpha - \{ 1 \}) \cdot (U_{\gamma_j} - \{ 1 \}) \bigr) \cap K$.
Therefore, we consider an element $v \in K$ of the form $v = u_\alpha(r)u_{\gamma_j}(s)$ with $r,s \in F^\times$. 
It remains to use the action of the torus $\mathscr{T}_A(F)$ to separate the two factors $U_\alpha$ and $U_{\gamma_j}$.
Again we compute for $v$ as above and $t = h_\alpha(\mu) h_\beta(\nu)$: 

\medskip

\centerline{$[t,v] = (t u_\alpha(r)u_{\gamma_j}(s) t^{-1}) \bigl(u_\alpha(r)u_{\gamma_j}(s)\bigr)^{-1} 
= u_{\alpha}\bigl( \alpha(t)r \bigr) u_{\gamma_j}\bigl( \gamma_j(t) s \bigr) u_{\gamma_j}(-s) u_\alpha(-r)$.} 

\medskip

\noindent
In view of the previous computation, and since $U_\alpha$ and $U_{\gamma_j}$ commute (this is where we use $m,n \geqslant 2$), this provides:

\medskip

\centerline{$[t,v] = u_{\alpha}\bigl( ( \mu^2 \nu^{-m}-1) r \bigr) u_{\gamma_k} \bigl( (\mu^{2a_j-nb_j} \nu^{-ma_j+2b_j}-1) s \bigr)$.} 

\medskip

\noindent
Now we can specialize our choice of multiplicative parameters $\mu$ and $\nu$.
For $\kappa \in F^\times$ we set $\mu = \kappa^m$ and $\nu = \kappa^2$; then for $t = h_\alpha(\kappa^m) h_\beta(\kappa^2)$ we obtain:

\medskip

\centerline{$[t,v] = u_{\gamma_j} \bigl( (\kappa^{m(2a_j-nb_j)} \kappa^{2(-ma_j+2b_j)}-1) s \bigr) = u_{\gamma_j} \bigl( (\kappa^{(4-mn)b_j}-1) s \bigr)$.} 

\medskip

\noindent
It remains to choose $\kappa \in F$ so that $\kappa^{(4-mn)b_j} \neq 1$ to conclude that $K \cap U_{\gamma_j} \neq \{ 1 \}$ and
$K \cap U_\alpha \neq \{ 1 \}$.

Similarly we can obtain $K \cap U_\beta \neq \{ 1 \}$. Therefore,
again using the action of $\mathscr{T}_A(F)$, we obtain
$U_\alpha, U_\beta \subset K$,
which finally shows that $K = \mathscr{G}_A(F)$.
\hfill$\square$

\medskip

\begin{remark}
\label{rk - simple}
Let us explain here why the same simplicity result over $F$ is easier when simplicity over finite fields is known.
Indeed let $\mathscr{G}_A$ be a simply connected Kac-Moody group for which simplicity is known over (sufficiently large) finite fields and let $K \triangleleft \mathscr{G}_A(F)$ be non-central.
Then, arguing as in the beginning of the above proof, we know that there exists $\ell \geqslant 1$ such that 
$K \cap \Bigl( \mathscr{G}_A({\bf F}_{q^\ell}) - Z\bigl(  \mathscr{G}_A({\bf F}_{q^\ell}) \bigr) \Bigr) \neq \varnothing$. 
Up to enlarging $\ell$, simplicity of $\mathscr{G}_A({\bf F}_{q^\ell})/ Z\bigl(  \mathscr{G}_A({\bf F}_{q^\ell}) \bigr)$ implies that $K$ contains the latter group, hence intersects non-trivially all the root groups, which finally implies that $K = \mathscr{G}_A(F)$.
\end{remark}

\bigskip

\subsection{Simplicity using simplicity}
\label{ss - heuristic}
For the sake of completeness, we conclude by explaining how simplicity for hyperbolic rank 2 Kac-Moody groups with non-trivial commutation relations for prenilpotent pairs can be proved.

\smallskip 

Recall that if $\Gamma$ is an infinite finitely generated group satisfying (NSP), then $\Gamma/Z(\Gamma)$ is called {\it just infinite}~in the sense that all its proper quotients are finite; this is, so to speak, half of simplicity (\ref{ss - NSP}).
Recall also that  for an infinite finitely generated group, the following implications are well-known: linearity $\Rightarrow$ residual finiteness $\Rightarrow$ non-simplicity (a group $\Gamma$ is said to be {\it residually finite}~if we have $\bigcap_{[\Gamma:\Delta]<\infty} \Delta = \{1\}$).
Here is a rough strategy to construct simple groups.
Let $\Gamma$ be an infinite group acting geometrically on a \cat-space.
Assume in addition that $\Gamma$ is both just infinite and {\it not}~residually finite.
Then the normal subgroup $\Gamma^\circ = \bigcap_{[\Gamma:\Delta]<\infty} \Delta$ is non-trivial, so it is a finite index subgroup since $\Gamma$ is just infinite.
In fact, more can be said: $\Gamma^\circ$ is a finite direct product of simple groups (all isomorphic to one another) \cite{Wilson}. 
It remains then to stand by the geometric situation (e.g. a suitable irreducibility of the geometric action) to be able to conclude that $\Gamma^\circ$ contains only one factor.

By (NSP), we know that a Kac-Moody lattice is just infinite (modulo center).
Therefore it is enough to show that a non-affine Kac-Moody lattice is non-residually finite, for instance because it contains a suitable non-residually finite subgroup. 
The latter subgroup can be given by some wreath product: {\it if $F$ is a finite non-abelian group, then $F \wr {\bf Z} = F^{({\bf Z})} \rtimes {\bf Z}$ is not residually finite}
\cite{Meskin}.
Using this, the following simplicity theorem can be proved \cite{CaReRk2}.

\begin{thm}
\label{thm - general}
Let $A = \left( \begin{array}{cc} \hfill 2  & -n \\ -1 & \hfill 2 \end{array} \right)$ be a generalized Cartan matrix of indefinite type (i.e. $n>4$),
and let $F$ be an algebraic closure of a finite field ${\bf F}_q$. Then, the corresponding simply connected Kac-Moody groups $\mathscr{G}_A({\bf F}_q)$ and $\mathscr{G}_A(F)$ are simple groups modulo their centers.
\end{thm}

\noindent {\it Reference}. This is \cite[Theorem 2]{CaReRk2}.
\hfill$\square$

\medskip

\noindent Summarizing all (including known) facts, and taking into account Remark \ref{rk - simple}, 
we obtain the following statement.

\begin{remark}
\label{rk - simple final}
Let $A$ be a generalized Cartan matrix of non-affine type, and let $\mathscr{G}_A$ be a Tits functor of type $A$.
Let $G$ be the elementary subgroup of $\mathscr{G}_A(F)$ over
the algebraic closure $F$ of a finite field ${\bf F}_p$ (that is, $G = [ \mathscr{G}_A(F) , \mathscr{G}_A(F) ]$); the group
$G$ is generated by all root subgroups. Then, $G$ is a simple group modulo its center whenever $A$ is indecomposable. 
\end{remark}

At last, it is natural to formulate the following question.

\begin{question}
Let $A = \left( \begin{array}{cc} \hfill 2  & -n \\ -m & \hfill 2 \end{array} \right)$ be a generalized Cartan matrix of indefinite type, i.e. $mn>4$.
Let $\mathscr{G}_A$ be the corresponding simply connected incomplete Kac-Moody group and let ${\bf F}_q$ be a finite field.
Assume that $m,n \geqslant 2$.
Is the finitely generated group $\mathscr{G}_A({\bf F}_q)/Z\bigl( \mathscr{G}_A({\bf F}_q) \bigr)$ simple?
\end{question}

Simplicity in this case would shortcut the proof of the present paper, but we think that providing a simplicity proof over $\overline{{\bf F}_q}$, using only the weakening of simplicity (NSP) over ${\bf F}_q$, has its own interest.
Note that the above question also applies to more exotic lattices of locally finite Moufang twin trees, as defined in \cite{AbRe}.
Some of these groups can be constructed with a trivial torus, which might be an obstruction to simplicity.

\end{document}